\documentclass[11pt]{amsart}

\usepackage{url}
\usepackage[vmargin=2.5cm, hmargin=2.5cm]{geometry}
\usepackage{kpfonts}
\usepackage{xcolor}

\title{\texttt{numericalsgps}, a \texttt{GAP} package for numerical semigroups}

\author{M. Delgado}
\address{CMUP, Departamento de Matematica, Faculdade de Ciencias, Universidade do Porto, Rua do Campo Alegre 687, 4169-007 Porto, Portugal}
\email{mdelgado@fc.up.pt}
\thanks{The authors were partially supported by CMUP (UID/MAT/00144/2013), which is funded by FCT (Portugal) with national (MEC) and European structural funds through the programs FEDER, under the partnership agreement PT2020.}

\author{P. A. Garc\'{\i}a-S\'anchez}
\address{Departamento de \'Algebra and CITIC-UGR, Universidad de Granada, E-18071 Granada, Espa\~{n}a}
\email{pedro@ugr.es}
\thanks{The second author is supported by the projects MTM2010-15595, FQM-343,  FQM-5849, NSF-1061366 and FEDER funds. Some functionalities have been tested in the supercomputer alhambra.ugr.es}

\makeatletter
\def\@setauthors{%
  \begingroup
  \def\thanks{\protect\thanks@warning}%
  \trivlist
  \centering\footnotesize \@topsep30\p@\relax
  \advance\@topsep by -\baselineskip
  \item\relax
  \author@andify\authors
  \def\\{\protect\linebreak}%
  \authors%
  \ifx\@empty\contribs
  \else
    ,\penalty-3 \space \@setcontribs
    \@closetoccontribs
  \fi
  \endtrivlist
  \endgroup
}
\def\@settitle{\begin{center}%
  \baselineskip14\p@\relax
    \bfseries
  \@title
  \end{center}%
}
\makeatother



\begin{document}

\begin{abstract}
The package \texttt{numericalsgps} performs computations with and for numerical semigroups. Recently also affine semigroups are admitted as objects for calculations. This manuscript is a survey of what the package does, and at the same time of the trending topics on numerical semigroups.
\end{abstract}

\maketitle

\section{Introduction}

The motivation for the implementation of \texttt{numericalsgps} was the lack of a package specific to make computations related to numerical semigroups. We had several functions implemented in distinct programming languages, each with its own interface, which made the communication between them very difficult. Thus we decided to unify all these procedures in a single package written in an appropriate language. The choice for the language was \texttt{GAP} (\cite{gap}), mainly because of the experience achieved by the first author. 

The first version of \texttt{numericalsgps} was released in 2005. Since then, the authors of this note have been adding new functionalities and replacing, when appropriate, algorithms with newer and faster ones. 
This makes the package to permanently reflect the state of the art in the area.

Also some algorithms have different implementations that are used taking into account either the information stored in the object to deal with, or what packages has the user installed/loaded in his \texttt{GAP} installation. 

The current version is 1.0 and can be found in \cite{numericalsgps}; the development version is available in \url{https://bitbucket.org/gap-system/numericalsgps}. The manual is over 90 pages long. 

The heart of this manuscript is Section~\ref{sec:Contents}, which consists of a brief description of the contents of the package. Aiming to make the paper self contained, we give definitions when necessary. The few examples given are simple illustrations that can guide the reader to produce his owns (assuming he has a working installation of GAP). We give pointers to the literature where one can find the implemented algorithms, which makes also the paper a kind of review of the computational procedures to deal with numerical and affine semigroups.

The paper ends with a reference to interactions with other commutative algebra packages. The use of external software frequently improves very much the execution time of the functions offered for affine semigroups.

\section{Contents}\label{sec:Contents}

We briefly describe in this section the contents of the package, using as a guideline the chapters of the manual.

\subsection{Introduction}

In the introduction of the manual, the basic definitions concerning numerical semigroups are given. The author interested in the topic can find all these definitions in \cite{ns}.

\subsection{Numerical semigroups}

A \emph{numerical semigroup} $S$ is a submonoid of the set of nonnegative integers $\mathbb N$ under addition, and such that $\mathbb N\setminus S$ has finitely many elements. This section describes several ways to define a numerical semigroup.

The elements in the set $\mathbb N\setminus S$ are usually called \emph{gaps}, and its cardinality is the \emph{genus} of $S$. We say that a gap $g$ is a \emph{fundamental gap} of $S$ if it is maximal in the set of gaps with respect to division, or in other words, $g\not\in S$, $2g\in S$ and $3g\in S$. Gaps and fundamental gaps fully determine the semigroup $S$, and so, they can be used to describe a numerical semigroup. Functions are provided to determine if a list of nonnegative integers is a list of gaps or fundamental gaps; and also procedures to define numerical semigroups by means of these lists.

Since $\mathbb N\setminus S$ is finite, the maximum of $\mathbb Z\setminus S$ exists and it is known as the \emph{Frobenius number} of $S$ (there is actually a huge number of papers dealing with the computation/bounds of the Frobenius number of numerical semigroups; see for instance \cite{alfonsin}). The \emph{conductor} of $S$ is just the Frobenius number of $S$ plus one, and has the property that it is the least nonnegative integer $c$ such that $c+\mathbb N\subseteq S$. We call the elements in $S$ less than or equal to the conductor the \emph{small elements} of $S$. Clearly, the semigroup $S$ is uniquely determined by its small elements. A procedure is implemented to check if a list of integers is the set of small elements of a numerical semigroup, and also a function to define a numerical semigroup if this is the case in terms of this list. 

If we take a closed interval $I=[a/b,c/d]$ with $a$, $b$, $c$ and $d$ positive integers such that $a/b<c/d$, then the set $ \bigcup_{k\in\mathbb N}(\mathbb N\cap kI)$ is a numerical semigroup (and coincides with the set of all numerators of rational elements in $I$). It can be shown that this class of semigroups coincides with the numerical semigroups that consist of nonnegative integer solutions to inequalities of the form $\alpha x\bmod \beta \le \gamma x$, which are known as \emph{proportionally modular} numerical semigroups. If $\gamma=1$, then they are simply called \emph{modular}. Hence we can also define a numerical semigroup in terms of the modular or proportionally modular inequality (giving a list with the parameters) or by an interval (providing its ends). Note that distinct intervals can yield the same numerical semigroup (and the same holds for proportionally modular inequalities). Membership to a numerical semigroup given by a proportionally modular inequality is trivial. Also specific fast algorithms exist for the computation of the Frobenius number if this is the case. For some kind of semigroups where testing being proportionally modular is fast, we perform this test and keep the inequality.

Another way to uniquely determine a numerical semigroup is by any of its Ap\'ery sets of its nonzero elements. Let $S$ be a numerical semigroup and let $n\in S\setminus\{0\}$. The \emph{Ap\'ery set} of $n$ in $S$ is the set $\{s\in S\mid s-n\not\in S\}$. This set has precisely $n$ elements, one for each congruence class modulo $n$. Once we know an Ap\'ery set, the cost of membership to $S$ is small, and also the Frobenius number and genus can be easily computed. Thus if the Ap\'ery set with respect to the least positive integer in $S$ (its \emph{multiplicity}) is computed, we store it as part of the object $S$. Many other invariants depend also on this specific Ap\'ery set as we will refer later. We provide a function to determine if a given list of integers is the Ap\'ery set of an element (the length of the list) in a numerical semigroup, and also to define a numerical semigroup by means of the Ap\'ery set.

Ap\'ery sets and proportionally modular inequalities can be seen as particular cases of periodic subadditive functions. We say that $f:\mathbb N\to \mathbb N$ is subadditive if $f(0)=0$ and $f(i+j)\le f(i)+f(j)$ for all $i,j\in \mathbb N$. Associated to $f$ we can define the semigroup of nonnegative integers $x$ such that $f(x)\le x$. This set is a numerical semigroup when $f$ is periodic (with positive period). We represent a periodic function by the values on the integers less than the period, and thus by a list of nonnegative integers. We give a function to test if a list corresponds to a subadditive function, and if so, a numerical semigroup can be defined by using this list as argument.

Let $A$ be a nonempty subset of $\mathbb N$. The monoid generated by $A$, denoted $\langle A\rangle$, is the set of all (finite) sums of elements of $A$. We say that $A$ \emph{generates} the numerical semigroup $S$ if $\langle A\rangle =S$. Observe that if this is the case, then the multiplicity of $S$ must be in $A$, and whenever two generators are congruent modulo the multiplicity, we do not need the largest one to generate the same semigroup. So we can always think of $A$ to be finite (since all its elements can be chosen to be incongruent modulo the multiplicity). Clearly, $S$ is uniquely determined by any of its systems of generators. Among these, there is only one minimal with respect to inclusion (actually also with respect to cardinality) which is $S^*\setminus (S^*+S^*)$, where $S^*=S\setminus\{0\}$. The cardinality of this set is known as the \emph{embedding dimension} of $S$. We give functions to define a numerical semigroup in terms of a generating set.

\begin{verbatim}
gap> s:=NumericalSemigroup("interval",71/5,153/8);
<Proportionally modular numerical semigroup satisfying 765x mod 10863 <= 197x >
gap> t:=NumericalSemigroup(15, 16, 17, 18, 19, 29, 43);
<Numerical semigroup with 7 generators>
gap> s=t;
true
\end{verbatim}

\subsection{Basic operations}

Among the basic operations of a numerical semigroup related to the contents of the preceding section, the package offers: computation of the multiplicity, generating system, minimal system of generators, small elements, gaps, embedding dimension, Ap\'ery sets, Frobenius number, conductor and fundamental gaps. 

Some functions have different methods depending on what is known about the semigroup. As an example, if the Ap\'ery set is known, the Frobenius number will be computed by using Selmer's formulas (see for instance \cite[Proposition 2.12]{ns}).

Given $S$ a numerical semigroup, we also give a procedure to list the first $n$ integers in $S$, with $n$ a positive integer. 

Associated to $S$ we can define the partial order relation $\le_S$ on $\mathbb Z$ as follows. We write 
\[a\le_S b \hbox{ if } b-a\in S.\] 
The set of maximal elements of $\mathbb Z\setminus S$ with respect to this order are known as \emph{pseudo-Frobenius} numbers (actually the Frobenius number is one of them), and their cardinality is the \emph{type} of $S$.  We provide functions to compute the pseudo-Frobenius numbers (that can be also obtained from the Ap\'ery sets) and the type of $S$.

Let $m$ be the multiplicity of $S$. Then the elements in the Ap\'ery set of $m$ in $S$ are $w_i=k_im+i$ for $i\in\{0,\ldots, m-1\}$ with $k_0=0$ and $(k_1,\ldots, k_{m-1})$ fulfilling a set of inequalities (\cite{london}). In this way a numerical semigroup with multiplicity $m$ corresponds with a point inside a polytope. We give a function that outputs the set of inequalities describing this polytope, and also to compute $(k_1,\ldots,k_{m-1})$, which are known as the Kunz coordinates of $S$.

An element $s\in S$ is a minimal generator if $S\setminus \{s\}$ is again a numerical semigroup. Hence the dual of this property could be an element $g\not \in S$ such that $S\cup\{g\}$ is also a numerical semigroup. These elements are known as \emph{special gaps}. We give a function to compute them, that can be used to compute oversemigroups of a given semigroup (Section \ref{constructing}).

\subsection{Presentations of a numerical semigroup} \label{presentations}
Let $S$ be a numerical semigroup minimally generated by $\{n_1,\ldots, n_e\}$. Then the monoid morphism $\varphi:\mathbb N^e\to S$, $\varphi(a_1,\ldots , a_e)=\sum_{i=1}^e a_i n_i$ is an epimorphism, known as the \emph{factorization homomorphism} of $S$. Consequently $\mathbb N^e/\ker\varphi$ is isomorphic to $S$, where $\ker\varphi=\{ (a,b)\in \mathbb N^e\times \mathbb N^e\mid \varphi (a)=\varphi(b)\}$. A \emph{presentation} of $S$ is a generating system of the congruence $\ker\varphi$. A \emph{minimal presentation} of $S$ is a minimal generating system of $\ker\varphi$ (again, no matter if you think about minimal with respect to inclusion or to cardinality; both concepts coincide for numerical semigroups; see \cite[Chapter 7]{ns}). 

Minimal presentations can be computed from graphs associated to elements in the numerical semigroup. Let $n$ be a nonzero element of $S$. We define the graph associated to $n$ as the graph with vertices the generators $n_i\in\{n_1,\ldots, n_e\}$ such that $n-n_i\in S$; and $n_in_j$ is an edge if $n-(n_i+n_j)\in S$. There is a function to compute the graph associated to $n$. A minimal presentation is constructed from those graphs that are not connected (there are finitely many of them and can be found by using, once more the Ap\'ery set of the multiplicity).  The elements having an associated non connected graph are called \emph{Betti elements} of $S$. A procedure to find the set of Betti elements of $S$ is given in the package; and also to find a minimal presentation of $S$.

Some numerical semigroups admit essentially a unique minimal presentation, in the sense that if $\sigma$ and $\tau$ are two minimal presentations (and thus have the same cardinality), whenever $(a,b)\in \sigma$, either $(a,b)\in \tau$ or $(b,a)\in \tau$ (that is, unique up to permutation of the pairs of the presentation). In particular, generic numerical semigroups have unique minimal presentations (\cite{b-gs-g}). The semigroup $S$ is \emph{generic} if every pair $(a,b)$ in a minimal presentation of $S$ has the property that $a-b$ has no zero coordinates. We give procedures to detect whether or not $S$ is uniquely presented or generic.

A straight generalization of the graph associated to $n\in S$ is the following: we can construct the simplicial complex of subsets $A$ of $\{n_1,\ldots, n_e\}$ such that $n-\sum_{a\in A}a \in S$. This set is known as the shaded set of $n$ in $S$ and has some nice properties associated to the generating function of $S$ (\cite{s-w}).

The congruence $\ker\varphi$ is also a submonoid of $\mathbb N^e\times \mathbb N^e$, which is generated by its nonzero minimal elements with respect to the usual partial ordering on $\mathbb N^e\times \mathbb N^e$ (\cite[Chapter 8]{fg}). If $(a,b)$ is one of this minimal generators, then $\varphi(a)=\varphi(b)\in S$ is called a \emph{primitive element} of $S$. These elements play an important role in factorization properties of $S$, and consequently we provide a function to compute them.

\begin{verbatim}
gap> s:=NumericalSemigroup(5,7,9);
<Numerical semigroup with 3 generators>
gap> MinimalPresentationOfNumericalSemigroup(s);
[ [ [ 0, 2, 0 ], [ 1, 0, 1 ] ], [ [ 4, 1, 0 ], [ 0, 0, 3 ] ], 
[ [ 5, 0, 0 ], [ 0, 1, 2 ] ] ]
\end{verbatim}

\subsection{Constructing numerical semigroups from others}\label{constructing}
We have already seen that adding a special gap to a numerical semigroup produces a new numerical semigroup, and the same holds if we remove a minimal generator. The intersection of two numerical semigroups also produces a numerical semigroup. Functions performing these tasks are provided in \texttt{numericalsgps}.

Let $p$ be a positive integer and $S$ be a numerical semigroup. The set $S/p=\{x\in \mathbb N\mid px \in S\}$ is again a numerical semigroup, called the \emph{quotient} of $S$ by $p$. A function is given to compute this new semigroup. 
A kind of inverse is the notion of \emph{multiple} of a numerical semigroup: given an integer $a>1$ and a numerical semigroup $S$, then $aS$ is a submonoid of $\mathbb N$, but it is not a numerical semigroup. If we add to this set all the integers greater than or equal to a given positive integer, say $b$, then we obtain a numerical semigroup: $aS\cup \{b,\to\}$. If we start from $\mathbb N$, and we repeat this operation several times, then we construct the set of what is known in the literature as \emph{inductive numerical semigroups} (see for instance \cite{f-gs} and the references therein) .

For a numerical semigroup $S$ the set of numerical semigroups $T$ with $S\subseteq T\subseteq \mathbb N$ is finite (the \emph{oversemigroups} of $S$), since the genus of $S$ is always finite by definition. We provide a function to compute the set of all oversemigroups of a given semigroup. Also there is a procedure to compute all numerical semigroups with given Frobenius number (this is done using the concept of fundamental gap as explained in \cite{fund-gap}) and another function to compute the set of all numerical semigroups with given genus $g$ (by constructing the tree of all numerical semigroups up to the level $g$).

\begin{verbatim}
gap> s:=NumericalSemigroup(5,7,9);
<Numerical semigroup with 3 generators>
gap> Length(OverSemigroupsNumericalSemigroup(s));
15
gap> Length(NumericalSemigroupsWithFrobeniusNumber(21));
1828
\end{verbatim}

We provide functions implementing the algorithms given in \cite{DGSRP15} to compute the set of all numerical semigroups having a given set as set of pseudo-Frobenius numbers.

\begin{verbatim}
gap> pf := [13,24,25];;                                
gap> NumericalSemigroupsWithPseudoFrobeniusNumbers(pf);
[  ]
gap> pf := [13,19,25];;                                
gap> NumericalSemigroupsWithPseudoFrobeniusNumbers(pf);
[ <Numerical semigroup>, <Numerical semigroup>, <Numerical semigroup> ]
\end{verbatim}

\subsection{Irreducible numerical semigroups}\label{irreducibles}

A numerical semigroup is \emph{irreducible} if it cannot be expressed as the intersection of two numerical semigroups properly containing it. This is equivalent to saying that it is maximal in the set of numerical semigroups with its same Frobenius number. Every numerical semigroup can be expressed (though not uniquely) as an intersection of irreducible numerical semigroups. We give a function to do this in \texttt{numericalsgps} (see \cite[Chapter 3]{ns} for a description of the algorithm). 

We also give a procedure to compute all irreducible numerical semigroups with given Frobenius number: the procedure is based in \cite{bl-r-irr}. This is actually equivalent to computing all irreducible numerical semigroups with given genus. This is due to the fact that if $f$ is the Frobenius number of an irreducible numerical semigroup, then either $g=(f+1)/2$ or $g=(f+2)/2$, depending on the parity of $f$. 

A numerical semigroup $S$ with Frobenius number $f$ is \emph{symmetric} if whenever $x\in\mathbb Z\setminus S$, $f-x\in S$.  The class of symmetric numerical semigroups coincides with that of irreducible numerical semigroups with odd Frobenius number. Irreducible numerical semigroups with even Frobenius number are called \emph{pseudo-symmetric}. We give tests to detect if a numerical semigroup is in any of these classes. 

A particular class of irreducible numerical semigroups is the set of numerical semigroups with the least possible number of relations in its minimal presentations. These semigroups are called \emph{complete intersections}, and it can be shown that every complete intersection numerical semigroup is either $\mathbb N$ or a gluing of two complete intersections (see for instance \cite[Chapter 8]{ns}). We say that $S=S_1+S_2$, with $S$ a numerical semigroup and  $S_1$ and $S_2$ submonoids of $\mathbb N$, is a \emph{gluing} of $S_1$ and $S_2$ if $\gcd(S_1)\gcd(S_2)\in S_1\cap S_2$ and $\gcd(S_1)\neq 1\neq \gcd(S_2)$. We give procedures to detect if a numerical semigroup can be expressed as a gluing of two of its submonoids, and if it is a complete intersection. 

We also implement the procedures presented in \cite{ci} to compute the set of all complete intersection numerical semigroups with fixed Frobenius number (equivalently fixed genus, since we are still dealing with irreducible numerical semigroups). We present procedures to detect if a numerical semigroup is \emph{free} (either $\mathbb N$ or a gluing of a free numerical semigroup with a copy of $\mathbb N$) and to calculate all free numerical semigroups with fixed Frobenius number. The same is done for \emph{telescopic} numerical semigroups (these are free numerical semigroups where the gluing is performed in the same order given by the generators) and numerical semigroups associated to irreducible planar curve singularities (a particular case of telescopic numerical semigroups; see \cite{ci} for more details).

A generalization of the concept of irreducible numerical semigroup is the following. We have seen that the genus $g$ of an irreducible numerical semigroup $S$ with Frobenius number $f$ is either $g=(f+1)/2$ if $f$ is odd (symmetric), or $g=(f+2)/2$ if $f$ is even (pseudo-symmetric). It turns out that the type of symmetric numerical semigroups is 1 and the type of pseudo-symmetric numerical semigroups is 2. So if $S$ is an irreducible numerical semigroup with genus $g$, Frobenius number $f$ and type $t$, then $g=(f+t)/2$. We say that a numerical semigroup $S$ is \emph{almost-symmetric} if its genus is one half of its Frobenius number plus its type. We give a function to test if a numerical semigroup is almost-symmetric and include the procedure presented in \cite{almost} to compute the set of almost symmetric numerical semigroups with fixed Frobenius number.

\begin{verbatim}
gap> s:=NumericalSemigroup(5,7,9);
<Numerical semigroup with 3 generators>
gap> DecomposeIntoIrreducibles(s);
[ <Numerical semigroup>, <Numerical semigroup> ]
gap> List(last,MinimalGeneratingSystem);
[ [ 5, 7, 8, 9 ], [ 5, 7, 9, 11 ] ]
gap> Length(TelescopicNumericalSemigroupsWithFrobeniusNumber(101)); 
86
gap> Length(AlmostSymmetricNumericalSemigroupsWithFrobeniusNumber(31));
1827
\end{verbatim}

\subsection{Ideals of numerical semigroups}\label{ideals}

A nonempty subset $I$ of $\mathbb Z$ is a \emph{relative ideal} of a numerical semigroup $S$ if $I+S\subseteq I$ and there exists $d\in \mathbb Z$ such that $d+I\subseteq S$ (the concept of relative ideal corresponds to that of fractional ideal in domains). Every relative ideal $I$ of $S$ can be expressed in the form $I=\{i_1,\ldots, i_n\}+S$ for some integers $i_j$. The set $\{i_1,\ldots, i_n\}$ is a \emph{generating set} of the ideal, and it is minimal if no proper subset generates the same ideal. 

\begin{verbatim}
gap> s:=NumericalSemigroup(3,4,5);
<Proportionally modular numerical semigroup satisfying 5x mod 15 <= 2x >
gap> 5+s;
<Ideal of numerical semigroup>
gap> [-1,2]+s;
<Ideal of numerical semigroup>
gap> MinimalGeneratingSystem(last);
[ -1 ]
\end{verbatim}

We provide functions for computing the small elements of an ideal (the definition is analogous to that in numerical semigroups), Ap\'ery sets (and tables; see \cite{cbjza13}), the ambient numerical semigroup, membership, and also some basic operations as addition, union, subtraction ($I-J=\{z\in \mathbb Z\mid z+J\subseteq I\}$), set difference, multiplication by an integer, translation by an integer, intersection, blow-up ($\bigcup_{n\in \mathbb N} nI-nI$) and $*$-closure with respect to a family of ideals (\cite{spi}). 

Numerical semigroups are ``local'' in the sense that there is a unique maximal ideal: the set of nonzero elements of the semigroup. Also there exists a \emph{canonical ideal}, which for a numerical semigroup $S$ with Frobenius number $f$ is defined as $\{z\in \mathbb Z\mid f-z\not\in S\}$ (see for instance \cite{bf06}).

The \emph{Hilbert function} associated to an ideal $I$ of a numerical semigroup $S$ is the function that maps every $n\in \mathbb N$ to $nI\setminus(n+1)I$. The \emph{reduction number} of $I$ is the least positive integer $n$ such that $\min(I)+nI=(n+1)I$. We give functions to compute the reduction number and the Hilbert function associated to an ideal. Also we give a procedure that computes the microinvariants of a numerical semigroup which are used to determine if the graded ring associated to the semigroup ring $K\llbracket S\rrbracket$ is Cohen-Macaulay (see \cite{bf06}). 

Finally we give a function to test if a numerical semigroup is a monomial semigroup ring following \cite{mi02}. A numerical semigroup $S$ is said to be \emph{monomial} if for any ring $R$ with $K\subseteq R\subseteq K\llbracket x\rrbracket$ and such that the algebraic closure of $R$ is $K\llbracket x\rrbracket$ ($K$ a field with characteristic zero) 
and $\mathrm v(R)=S$, we have that $R$ is a semigroup ring. Here, $\mathrm v$ denotes the usual valuation.

\subsection{Numerical semigroups with maximal embedding dimension}

Recall that the embedding dimension of a numerical semigroup is the cardinality of its unique minimal generating system. Clearly, two minimal generators cannot be congruent modulo the multiplicity of the semigroup (the least positive integer in the semigroup). As a consequence, the embedding dimension is at most the multiplicity of the semigroup. Thus we say that a numerical semigroup $S$ has \emph{maximal embedding dimension} if its embedding dimension equals its multiplicity. 

The set of maximal embedding dimension numerical semigroups with fixed multiplicity, say $m$, is closed under intersection, and also if $S\neq \{0\}\cup(m+\mathbb N)$, then $S\cup\{f\}$ is also a maximal embedding numerical semigroup, with $f$ the Frobenius number of $S$. This in particular implies that if we are given a numerical semigroup that is not of maximal embedding dimension, we can consider the set of all maximal embedding dimension numerical semigroups with its same multiplicity containing it, and then the intersection of all of them, obtaining in this way the maximal embedding dimension closure of the given semigroup. Following this idea one can define the concept of minimal generators with respect to this class: the elements in the semigroup so that the closure of them yields the given semigroup. These elements are precisely the elements $x$ in a maximal embedding dimension numerical semigroup $S$ (together with the multiplicity) such that $S\setminus\{x\}$ is a maximal embedding dimension numerical semigroup. 
\begin{verbatim}
gap> s:=NumericalSemigroup(3,5,7);;
gap> MinimalMEDGeneratingSystemOfMEDNumericalSemigroup(s);
[ 3, 5 ]
\end{verbatim}
We also give functions to compute the maximal embedding dimension closure of an arbitrary numerical semigroup.

If $S$ is a numerical semigroup with multiplicity $m$, then $S$ has maximal embedding dimension if and only if for every $x,y\in S\setminus\{0\}$, $x+y-m\in S$. A natural generalization of this pattern is the following. We say that a numerical semigroup $S$ is Arf if for any $x,y,z\in S$ with $x\ge y\ge z$, then $x+y-z\in S$. Clearly, every Arf numerical semigroup has maximal embedding dimension. Also, the class of Arf numerical semigroups is closed under finite intersections and the adjoin of the Frobenius number (of course if we are considering semigroups other than $\mathbb N$). Thus the class of Arf numerical semigroups is a Frobenius variety (\cite[Chapter 6]{ns}). Again, it makes sense to talk about minimal generators with respect to this class, and also about the Arf closure of a given numerical semigroup (the intersection of all Arf numerical semigroups containing it). We give functions computing both things: Arf minimal generating sets and Arf closures. Also we provide a method to detect if a numerical semigroup is Arf, and a procedure that calculates the set of all Arf numerical semigroups with given Frobenius number.

Finally, we consider in this section the class of saturated numerical semigroups, which turns out to be again a Frobenius variety (closed under intersections and the adjoint of the Frobenius number). A numerical semigroup is \emph{saturated} if for every $s,s_1,\ldots, s_r\in S$ with $s_i\le s$ for all $i$ and every $z_1,\ldots, z_r\in\mathbb Z$ such that $z_1s_1+\cdots +z_rs_r\ge 0$ one gets $s+z_1s_1+\cdots +z_rs_r\in S$. We provide for saturated semigroups the analogous functions that we described in the above paragraph for Arf semigroups.

\subsection{Nonunique invariants for factorizations in numerical semigroups}

Let $S$ be a numerical semigroup minimally generated by $\{n_1,\ldots, n_e\}$. Recall that we defined a monoid epimorphism in Section \ref{presentations}, $\varphi:\mathbb N^{e}\to S$, $\varphi(a_1,\ldots, a_e)=a_1n_1+\cdots +a_en_e$. Observe that for $s\in S$, $\mathsf Z(s)=\varphi^{-1}(s)$ collects the different expressions of $s$ in terms of the generators of $S$. Thus we say that $\mathsf{Z}(s)$ is the set of \emph{factorizations} of $s$. The cardinality of $\mathsf{Z}(s)$ is usually known as the \emph{denumerant} of $s$. We use \texttt{RestritcedPartitions} to compute the set factorizations of $s$. 

The number of connected components of the graph associated to $s\in S$ (Section \ref{presentations}) coincides with the number of connected components of the graph with vertices given by $\mathsf Z(s)$ and $zz'$ is an edge provided that $z\cdot z'\neq 0$. 

Given $z=(z_1,\ldots,z_e)$ a factorization of $s\in S$, we write $|z|$ to denote the \emph{length} of $z$, $|z|=z_1+\cdots +z_e$. The \emph{maximal denumerant} of $s$ is the number of factorizations of $s$ with maximal length. Even though the denumerant is not bounded while $s$ increases in $S$, the maximal denumerant is finite and can be effectively computed (\cite{max-den}). We include this algorithm in the package as well as tests for supersymmetry and additiveness (see \cite{max-den} for details).
	
Let $A$ be the Ap\'ery set of $n_e$ in $S$. A subset $L$ of $\mathbb N^{e-1}$ is an \textsf{L}shape associated to $S$ if (1) $L\subset \mathsf Z(A)$ (the set of factorizations of the elements in $A$), (2) for every $a\in A$, $\#(L\cap \mathsf Z(a))=1$, and (3) for every $l\in L$, if $l'\in \mathbb N^{e-1}$ is such that $l'\le l$, then $l'\in L$. These sets give information on the factorizations on numerical semigroups (\cite{lformas-fact}), and this is why we have included a procedure to compute them.

The set of lengths of factorizations of $s$ is always finite (due to Dickson's lemma) and consequently we can write it as $\{l_1<\cdots <l_t\}$. The set $\{l_2-l_1,l_3-l_2,\ldots, l_t-l_{t-1}\}$ is the \emph{Delta set} associated to $s$. The Delta set of $S$ is the union of all the Delta sets of $s$. This set is finite, and its maximum is achieved in one of the Betti elements of $S$ (\cite{deltas}).

The \emph{elasticity} of an element $s$ is the ratio between the maximal and minimal lengths of factorizations of $s$. It was introduced to measure how far is a domain from being half-factorial (all factorizations of all the elements have the same length). No numerical semigroup other than $\mathbb N$ is half-factorial, which is a unique factorization monoid. We give a procedure to compute this invariant. 

Given $z=(z_1,\ldots, z_e), z'=(z_1',\ldots, z_e')\in \mathsf Z(s)$, we denote by $z\wedge z'=(\min(z_1,z_1'), \ldots, \min(z_e,z_e'))$, which corresponds to the ``common part'' of these factorizations. The \emph{distance} between $z$ and $z'$ is the $\mathrm d(z,z')=\max(|z-z\wedge z'|, |z'-z\wedge z'|)$. The \emph{catenary degree} of $s$ is the least positive integer such that for any two factorizations of $s$, there exists a chain of factorizations such that the distance between two consecutive factorizations is bounded by this integer. The catenary degree of $S$ is defined as the supremum of the catenary degrees of its elements. This supremum is reached in one of its Betty elements (\cite{cat-tame}). We give procedures to compute the catenary degree of a set of factorizations and of a numerical semigroup. Also other variants of catenary degrees are included: adjacent, homogeneous, equal or monotone catenary degree (see \cite{g-hk, hom}). For the homogenization of a numerical semigroup we offer a series of auxiliary functions.

The \emph{tame degree} of $s\in S$ is the least positive integer $t$ such that for every factorization $z$ of $s$ and every integer $i\in\{1,\ldots, e\}$ such that $s-n_i\in S$, there exists another factorization $z'$ of $s$ with nonzero $i$th coordinate and such that the distance to $z$ is less than or equal to $t$ (there exists a factorization in which $n_i$ is involved at a distance at most $t$). The tame degree of the semigroup $S$ is the supremum of all the tame degrees of its elements, and it is reached in one of its primitive elements (also in an element with associated noncomplete graph). We give functions to compute the tame degree of a set of factorizations and that of the semigroup.

Recall that associated to the numerical semigroup $S$, we can define the partial order on $\mathbb Z$, $a\le_S b$ if $b-a\in S$. Thus $(\mathbb Z,\le_S)$ is a poset, and one can define the M\"obius function associated to it. We implement the procedure presented in \cite{mob}.

The last invariant we give procedures to compute is the $\omega$-primality, which determines how far is an element from being prime. The \emph{$\omega$-primality} of $s\in S$ is the least positive integer $\omega$ such that whenever $s \le_S \sum_{a\in A} a$ with $A\subseteq S$ finite, there exists $\Omega\subseteq A$ with $\#\Omega\le \omega$ such that $s\le_S \sum_{a\in \Omega}a$. Clearly, if the omega primality is one, then the element is prime, if we look at $\le_S$ as a division.
The $\omega$-primality of the semigroup is the maximum of the $\omega$-primalities of its minimal generators. Initially we used the algorithm presented in \cite{b-gs-g}. Now we use a faster procedure implemented by C. O'Neill (see Section \ref{contrib}).

\begin{verbatim}
gap> l:=FactorizationsIntegerWRTList(100,[10,11,13,15]); 
[ [ 10, 0, 0, 0 ], [ 1, 7, 1, 0 ], [ 3, 4, 2, 0 ], [ 5, 1, 3, 0 ], 
[ 0, 2, 6, 0 ], [ 3, 5, 0, 1 ], [ 5, 2, 1, 1 ], [ 0, 3, 4, 1 ], 
[ 2, 0, 5, 1 ], [ 7, 0, 0, 2 ], [ 0, 4, 2, 2 ], [ 2, 1, 3, 2 ], 
[ 0, 5, 0, 3 ], [ 2, 2, 1, 3 ], [ 4, 0, 0, 4 ], [ 1, 0, 0, 6 ] ]
gap> TameDegreeOfSetOfFactorizations(l);
5
gap> CatenaryDegreeOfSetOfFactorizations(l);
3
\end{verbatim}

\subsection{Polynomials, formal series and numerical semigroups}

Let $S$ be a numerical semigroup. The \emph{Hilbert series} (not to be confused with the Hilbert function in Section \ref{ideals}) is the formal series $\mathrm H_S(x)= \sum_{s\in S} x^s$. Clearly $\sum_{n\in \mathbb N} x^n=1/(1-x)=\sum_{s\in\mathbb N\setminus S}x^s+\mathrm H_S(x)$. Hence $\mathrm P_S(s)=1+(x-1)\sum_{s\in \mathbb N\setminus S} x^s=(1-x)\mathrm H_S(x)$ is a polynomial, which we call the \emph{polynomial associated} to $S$ (see \cite{moree}). We provide functions to compute both the polynomial and Hilbert series of a numerical semigroup.

It turns out that when $S$ is a complete intersection, the polynomial associated to $S$ has all its roots in the unit circumference (and zero is not a root, which by Kronecker's lemma means that all the roots are in the unit circle, or equivalently, it is a product of cyclotomic polynomials). We give functions to determine if a monic polynomial with integer coefficients has all its roots in the unit circle, and to do this we need two auxiliary implementations: that of being cyclotomic and the computation of the Graeffe polynomial (see \cite{c-gs-m} for details). A numerical semigroup is said to be \emph{cyclotomic} if its associated polynomial has all its roots in the unit circle. 

Symmetry (see Section \ref{irreducibles}) can also be characterized in terms of the associated polynomial: a numerical semigroup is symmetric if and only if its associated polynomial is self-reciprocal (a palindrome if we look at the coefficients).

Let $K$ be an algebraically closed field. And let $f\in K[x,y]$ represent an irreducible curve with one place at infinity. Take $g\in K[x,y]$ and set $\mathrm{int}(f,g)=\dim_K (K[x,y])/(f,g)$. Then the set $\{\mathrm{int}(f,g)\mid g\not \in (f)\}$ is a numerical semigroup. We give a procedure to implement it (see \cite{a-gs}). This kind of semigroups are generated by what is called a $\delta$-sequence. There is a function to compute all $\delta$-sequences with fixed Frobenius number (equivalently genus since these semigroups are complete intersections and thus symmetric). Also associated to any $\delta$-sequence there is a ``canonical'' planar curve, and we offer a method to compute it. 

Let $F$ be a set of polynomials. Then the set of values (respectively degrees) of the series (respectively polynomials) in the algebra $K\llbracket F\rrbracket$ (respectively $K[F]$) is a submonoid of $\mathbb N$. Under certain conditions it is a numerical semigroup, and we provide functions to compute it. Also to determine a basis of the algebra $K\llbracket F\rrbracket$ (or $K[F]$) such that the values (or degrees) minimally generate the semigroup of values of this algebra (see \cite{a-gs-m}).

\begin{verbatim}
gap> t:=Indeterminate(Rationals,"t");;
gap> l:=[t^4,t^6+t^7,t^13];
[ t^4, t^7+t^6, t^13 ]
gap> SemigroupOfValuesOfCurve_Local(l);
<Numerical semigroup with 4 generators>
gap> MinimalGeneratingSystem(last);
[ 4, 6, 13, 15 ]
gap> SemigroupOfValuesOfCurve_Local(l,"basis");    
[ t^4, t^7+t^6, t^13, t^15 ]
\end{verbatim}

\subsection{Affine semigroups}

An \emph{affine semigroup} is a finitely generated submonoid of $\mathbb N^n$ for some positive integer $n$. In the package, affine semigroups can be defined by means of generators, as the set of elements in the positive orthant of a subgroup of $\mathbb Z^n$ (full semigroups) or as the set of elements in the positive orthant of a cone (normal semigroups). Our intention is to provide as many functions as possible for affine semigroups as we offer for numerical semigroups. Along this line, we present methods for membership, computing minimal presentations, determine gluings, Betti and primitive elements, and the whole series of procedures for nonunique factorization invariants (an overview of the existing methods for the calculation of these invariants can be found in \cite{overview_non_unique}).  New procedures are now under development based in Hilbert functions and binomial ideals (\cite{chris}).

As an example, let us do some computations with $G\cap\mathbb N^3$, where $G$ is the subgroup of $\mathbb Z^3$ with defining equations $x+y\equiv 0\bmod 2$ and $x+z\equiv 0\bmod 2$ (this is actually the block monoid associated to $\mathbb Z_2^3$; see \cite{g-hk} for the definition of block monoid).
\begin{verbatim}
gap> a:=AffineSemigroup("equations",[[[1,1,0],[0,1,1]],[2,2]]);
<Affine semigroup>
gap> GeneratorsOfAffineSemigroup(a);
[ [ 0, 0, 2 ], [ 0, 2, 0 ], [ 2, 0, 0 ], [ 1, 1, 1 ] ]
gap> OmegaPrimalityOfAffineSemigroup(a);
3
gap> BettiElementsOfAffineSemigroup(a);
[ [ 2, 2, 2 ] ]
\end{verbatim}

\subsection{Random}\label{random}
Based on the the methods provided by GAP to create ``random'' objects, we provide some functions for ``random'' affine numerical semigroups. These are particularly useful to produce examples. Furthermore, they are extensively used each time new algorithms are implemented and tests need to be made. 
\begin{verbatim}
gap> l:=List([1..20], _->RandomNumericalSemigroup(5,200));;
gap> ls:=Filtered(l, s-> 1+FrobeniusNumber(s)=GenusOfNumericalSemigroup(s)*2);;
gap> List(ls,MinimalGeneratingSystem);
[ [ 8, 103 ], [ 25, 109 ], [ 35, 57, 125 ], [ 3, 52 ], [ 15, 170, 178 ], 
  [ 3, 145 ], [ 21, 68, 153 ] ]
\end{verbatim}

\subsection{Contributions}\label{contrib}

There is a special section devoted to contributions. So far we are happy to count with functions implemented by A. Sammartano and C. O'Neill (apart from those co-implementations with J. I. Garc\'\i a-Garc\'\i a and A. S\'anchez-R.-Navarro). 

The functions implemented by Sammartano are mainly focused on deriving properties of the semigroup algebra $k[[S]]$ and its associated graded algebra from properties of the numerical semigroup $S$. He offers procedures to determine purity and $M$-purity of $S$ (\cite{br}), Buchbaum (\cite{d-m-m}), Gorenstein (\cite{d-m-s}) and complete intersection (\cite{d-m-s-13b}) property for the graded algebra; some special shapes of the Ap\'ery sets ($\alpha$, $\beta$ and $\gamma$-rectangular, see \cite{d-m-s-13a}); and the type sequence of a numerical semigroup (\cite{b-d-f}).

O'Neill on his side offers methods dealing with non unique factorization invariants: factorizations, $\omega$-primality and Delta sets for a list of elements in a numerical semigroup, Delta sets for the whole semigroup, and periodicity for the Delta sets (\cite{b-on-p}).

\section{Interaction with other packages}

Since the first release of the package many other packages have come into scene (some still under development). We have tried to take advantage of this. Dealing with affine and numerical semigroups translates in many cases to computing nonnegative integer solutions of linear Diophantine equations or Gr\"obner basis calculations of binomial ideals. Hence the interaction with \texttt{singular} (\cite{Singular}), \texttt{Normaliz} (\cite{normaliz}) and \texttt{4ti2} (\cite{4ti2}) was a step forward for us. For \texttt{singular} there are several options to consider: \cite{singular-gap}, \cite{GradedModules} and SingularInteface \url{https://github.com/gap-system/SingularInterface}. As for \texttt{4ti2}, we can use \cite{4ti2Interface} and \cite{4ti2gap}, which is under development. Finally there is an interface for \texttt{Normaliz} that can be found in \url{https://github.com/fingolfin/NormalizInterface}. We have implemented different methods for each procedure depending on which of the above packages the user has loaded/installed.

%
%


\begin{thebibliography}{10}

\bibitem{4ti2} 4ti2 team, 4ti2--a software package for algebraic, geometric and combinatorial problems on linear spaces, available at \url{www.4ti2.de}.

\bibitem{lformas-fact} F. Aguil\'o-Gost, P. A. Garc\'\i a-S\'anchez, Factoring in embedding dimension three numerical semigroups, Electron. J. Comb. \textbf{17} (2010), \#R138.

\bibitem{ci} A. Assi, P. A. Garc\'{\i}a-S\'{a}nchez, Constructing the set of complete intersection numerical semigroups with a given Frobenius number, Applicable Algebra in Engineering, Communication and Computing \textbf{24}, (2013), 133--148.

\bibitem{a-gs} A. Assi and P. A. Garc\'\i a-S\'anchez, On curves with one place at infinity, arXiv:1407.0490, 2014.

\bibitem{a-gs-m} A. Assi, P. A. Garc\'\i a-S\'anchez, and V. Micale. Bases of subalgebras of $k\llbracket x\rrbracket$ and $k[x]$, arXiv, 1412.4089, 2014.

\bibitem{GradedModules}
M.~Barakat, S.~Gutsche, S.~Jambor, M.~Lange-Hegermann, A.~Lorenz, and
  O.~Motsak.
\newblock {GradedModules}, a homalg based package for the abelian category of
  finitely presented graded modules over computable graded rings, {V}ersion
  2014.09.17.
\newblock \url{http://homalg.math.rwth-aachen.de/~barakat/homalg-project/GradedModules},  Sep 2014.
\newblock GAP package.

\bibitem{b-d-f} V. Barucci, D. D. Dobbs, M. Fontana, Maximality properties in numerical semigroups and applications to one-dimensional analytically irreducible local domains,  Memoirs of the American Mathematical Society \textbf{598}, 1997.


\bibitem{bf06} V. Barucci, R. Fr\"oberg, Associated graded rings of one-dimensional analytically irreducible rings, J. Algebra \textbf{304} (2006), 349--358.

\bibitem{b-on-p} T. Barron, C. O'Neill, R. Pelayo, On the computation of delta sets and $\omega$-primality in numerical monoids, preprint, 2014.


\bibitem{b-gs-g} V. Blanco, P. A. Garc\'{\i}a-S\'{a}nchez, A. Geroldinger, Semigroup-theoretical characterizations of arithmetical invariants with applications to numerical monoids and Krull monoids, Illinois J. Math. \textbf{55} (2011), 1385--1414.

\bibitem{bl-r-irr} V. Blanco, J.C. Rosales, The tree of irreducible numerical semigroups with fixed Frobenius number, Forum Math.  \textbf{25} (2013), 1249--1261.

\bibitem{normaliz} W. Bruns, B. Ichim, T. R\"omer, and S\"oger C. Normaliz. algorithms for rational cones and affine monoids, \url{http://www.math.uos.de/normaliz}, 2014. 

\bibitem{br} L. Bryant, Goto numbers of a numerical semigroup ring and the Gorensteiness of associated graded rings, Comm. Algebra \textbf{38} (2010), 2092--2128.

\bibitem{max-den} L. Bryant, J. Hamblin, The maximal denumerant of a numerical semigroup, Semigroup Forum \textbf{86} (2013), 571--582.

\bibitem{mob} J. Chappelon, J. L. Ram\'\i rez Alfons\'\i n, On the M\"obius function of the locally finite poset associated with a numerical semigroup, Semigroup Forum \textbf{87} (2013), 313--330.

\bibitem{cat-tame} S. T. Chapman, P. A. Garc\'\i a-S\'anchez, D. Llena, V. Ponomarenko, J. C. Rosales, The catenary and tame degree in finitely generated commutative cancellative monoids, Manuscripta Math. \textbf{120} (2006), 253--264.


\bibitem{deltas} S. T. Chapman, P. A. Garc\'{\i}a-S\'{a}nchez, D. Llena, A. Malyshev, D. Steinberg, On the Delta set and the Betti elements of a BF-monoid, Arab J Math \textbf{1} (2012), 53--61.

\bibitem{c-gs-m} A. Ciolan, P.A. Garc\'\i a-S\'anchez, Cyclotomic numerical semigroups,  Max-Planck-Institut f\"ur Mathematik Preprint Series 2014 (64); arXiv:1409.5614, 2014.

\bibitem{singular-gap} M. Costantini and W. de Graaf, Gap package singular; the gap interface to singular,
\url{http://gap-system.org/Packages/singular.html}, 2012.


\bibitem{cbjza13} T. Cortadellas Ben\'{\i}tez, R. Jafari, S. Zarzuela Armengou. On the Ap\'ery sets of monomial curves. Semigroup Forum \textbf{86} (2013), 289--320.

\bibitem{d-m-m} M. D'Anna, M. Mezzasalma, and Micale V. On the buchsbaumness of the associated graded ring of a one-dimensional local ring. Comm. Algebra \textbf{37} (2009), 1594--1603.

\bibitem{d-m-s} M. D'Anna, V. Micale, and A. Sammartano. On the associated graded ring of a semigroup ring. J. Commut. Algebra \textbf{3} (2011), 147--168.

\bibitem{d-m-s-13a} M. D'Anna, V. Micale, A. Sammartano, Classes of complete intersection numerical semigroups, Semigroup Forum \textbf{88} (2014), 453--467.

\bibitem{d-m-s-13b} M. D'Anna, V. Micale, A. Sammartano, When the associated graded ring of a semigroup ring is complete intersection, J. Pure Appl. Algebra \textbf{217} (2013), 1007--1017.

\bibitem{Singular} W. Decker, G.-M. Greuel, G. Pfister, H. Sch\"onemann, SINGULAR 3-1-6 -- A computer algebra system for polynomial computations, \url{http://www.singular.uni-kl.de}, 2012.

\bibitem{numericalsgps} M. Delgado, P. A. Garc\'{\i}a-S\'{a}nchez, J. Morais,
\lq\lq NumericalSgps\rq\rq, A GAP package for numerical semigroups,
available via \url{http://www.gap-system.org}.

\bibitem{DGSRP15}
M.~{Delgado}, P.~A. {Garc{\'{\i}}a-S{\'a}nchez} and A.~M. {Robles-P{\'e}rez}, Numerical semigroups with a given set of pseudo-Frobenius numbers, arXiv:1505.08111, May 2015.

\bibitem{f-gs} J. I. Farr\'an and P. A. Garc\'{\i}a-S\'{a}nchez, The second Feng-Rao number for codes coming from inductive semigroups, arXiv:1505.01395.

\bibitem{gap}
  The GAP~Group, GAP -- Groups, Algorithms, and Programming, 
  Version 4.7.5,  
  2014,
  \url{http://www.gap-system.org}.

\bibitem{overview_non_unique} P. A. Garc\'\i a-S\'anchez, An overview of the computational aspects of nonunique factorization invariants, arxiv.org:1504.07424, April 2015.

\bibitem{4ti2gap} P.A. Garc\'{\i}a-S\'anchez and A. S\'anchez-R.-Navarro, 4ti2gap, GAP wrapper for 4ti2, \url{https://bitbucket.org/gap-system/4ti2gap}.

\bibitem{hom} P. A. Garc\'\i a-S\'anchez, I. Ojeda, A. S\'anchez-R.-Navarro , Factorization invariants in half-factorial affine semigroups, Internat. J. Algebra Comput. \textbf{23} (2013), 111--122.

\bibitem{g-hk} A. Geroldinger and F. Halter-Koch Non-Unique Factorizations, Chapman \& Hall/CRC, Boca Raton, FL, 2006.

\bibitem{4ti2Interface} S. Gutsche, 4ti2interface, a link to 4ti2, \url{http://www.gap-system.org/Packages/4ti2interface.html}, 2013.


\bibitem{mi02} V. Micale, On monomial semigroups, Communications in Algebra \textbf{30} (2002), 4687 -- 4698.

\bibitem{moree} P. Moree, Numerical semigroups, cyclotomic polynomials and Bernoulli numbers, Amer. Math. Monthly \textbf{121} (2014), 890--902.

\bibitem{chris} C. O'Neill, On factorization invariants and Hilbert functions, arXiv:1503.08351.

\bibitem{alfonsin} J. L. Ram\'{\i}rez Alfons\'{\i}n, The Diophantine Frobenius
Problem, Oxford University Press, 2005.

\bibitem{fg} {J. C.} Rosales and {P. A.} Garc\'{\i}a-S\'anchez, Finitely generated commutative monoids, Nova Science Publishers, Inc., New York, 1999.

\bibitem{ns} {J. C.} Rosales and {P. A.} Garc\'{\i}a-S\'anchez,  
     \newblock Numerical Semigroups, 
     \newblock Developments in Mathematics \textbf{20}, Springer, New York, 2009.

\bibitem{almost} J. C. Rosales, P. A. Garc\'{\i}a-S\'{a}nchez, Constructing almost symmetric numerical semigroups from almost irreducible numerical semigroups, Comm. Algebra \textbf{42} (2014), 1362--1367.

\bibitem{london} J. C. Rosales, P. A. Garc\'{\i}a-S\'anchez, J. I. Garc\'{\i}a-Garc\'{\i}a and M. B. Branco, Systems of inequalities and numerical semigroups, J. Lond. Math. Soc. \textbf{65}(3) (2002), 611--623.

\bibitem{fund-gap} J. C. Rosales, P. A.  Garc\'{\i}a-S\'anchez, J. I. Garc\'{\i}a-Garc\'{\i}a, and J. A. Jim\'enez-Madrid, Fundamental gaps in numerical semigroups, J. Pure Appl. Algebra \textbf{189} (2004), 301--313.

\bibitem{spi} D. Spirito, Star operations on numerical semigroups, Comm. Algebra, to appear, 2014.

\bibitem{s-w} L.A. Sz\'ekely and N.C. Wormald, Generating functions for the Frobenius problem with 2 and 3 generators, Math. Chronicle \textbf{15} (1986), 49--57.

\end{thebibliography}
\end{document}